\documentclass[12pt]{article}
\usepackage{amssymb,amsmath,latexsym}

\newtheorem{theorem}{Theorem}[section]
\newtheorem{proposition}[theorem]{Proposition}
\newtheorem{corollary}[theorem]{Corollary}

\newtheorem{example}[theorem]{Example}
\newtheorem{remark}[theorem]{Remark}

\newtheorem{definition}[theorem]{Definition}
\newtheorem{problem}[theorem]{Problem}
\newenvironment{proof}{\noindent{\sc Proof.}}{\hfill\qed}
\newcommand{\qed}{\quad\lower0.05cm\hbox{$\Box$}}

\newcommand{\D}{\mathcal{D}}
\newcommand{\N}{\mathcal{N}}
\newcommand{\Na}{\mathbb{N}}
\newcommand{\R}{{\mathbb R}}

\newcommand{\Law}{\mathcal{L}}
\newcommand{\M}{{\widetilde M}}
\newcommand{\VPN}{(V,\nu, \tau,\tau^*)}

\newcommand{\norm}[1]{\mid\hspace{-.06cm}\mid \hspace{-.1cm} #1 \hspace{-.1cm}\mid\hspace{-.06cm}\mid}
\newcommand{\abs}[1]{\mid \hspace{-.1cm} #1 \hspace{-.1cm}\mid}
\newcommand{\h}{\hat \;}

\title{Boundedness in generalized \v Serstnev PN spaces}

\author{
Bernardo Lafuerza-Guill\'{e}n\\ {\small Departamento de
Estad\'{\i}stica y Matem\'{a}tica Aplicada, Universidad de
Almer\'\i a,}\\ {\small E--04120 Almer\'\i a, Spain, e-mail: {\tt
blafuerz@ual.es}}\\
\and Jos\'{e} L. Rodr\'{\i}guez \\ {\small \'{A}rea de
Geometr\'{\i}a y Topolog\'{\i}a, Universidad de Almer\'\i a,}\\
{\small E--04120 Almer\'\i a, Spain, e-mail: {\tt
jlrodri@ual.es}}}
\date{\empty}

\begin{document}
\maketitle
\begin{abstract}
The motivation of this paper is a suggestion by H\"ole of
comparing the notions of $\D$-boundedness and boundedness in
Probabilistic Normed spaces (briefly PN spaces), with non
necessarily continuous triangle functions. Such spaces are here
called ``pre-PN spaces''. Some results on \v Serstnev spaces due
to B. Lafuerza, J. A. Rodr{\'\i}guez, and C. Sempi, are here extended
to {\it generalized \v Serstnev spaces} (these are pre-PN spaces
satisfying a more general \v Serstnev condition).

We also prove some facts on PN spaces (with continuous triangle
functions). First, a connection between fuzzy normed spaces
defined by Felbin and certain \v Serstnev PN spaces is
established. We further observe that topological vector PN spaces
are $F$-normable and paranormable, and also that locally convex
topological vector PN spaces are bornological. This last fact
allows to describe continuous linear operators between certain
generalized \v Serstnev spaces in terms of bounded subsets.

\medskip
\noindent {\it 2000  Mathematics Subject Classification:} 54E70,
46S50.
\end{abstract}

\section{Motivation and main results}

After the work done by Lafuerza, Rodr{\'\i}guez, and Sempi in
\cite{LRS3}, H\"ohle suggested the following:
\begin{problem}
\label{Hohle-problem} {\rm Compare $\D$-boundedness and
boundedness in the PN spaces $(V,\nu,\tau, \tau^*)$, with $\tau$
and $\tau^*$ non necessarily continuous.}
\end{problem}
Such spaces are here called {\it pre-PN spaces}  and include the
PN spaces introduced in \cite{alsina1} (where the triangle
functions $\tau$ and $\tau^*$ are assumed to be continuous). Every
PN space endows a topology, usually called the {\it strong
topology}, which is always metrizable (thanks to the continuity of
$\tau$ and $\tau^*$).
The probabilistic norm $\nu$ is a continuous map from $V$ to the
space $\Delta^+$ of distance distribution functions (the last
endowed with the Levy-Sybley metric). Recall that a base of strong
neighborhoods at $q\in V$ is $\{N_q(t)\}_{t>0}$ where
$$
N_q (t) :=\{p: \nu_{p-q}(t)> 1- t\} = q + N_\theta(t),
$$
$\theta$ being the origin of $V$.

In the case of pre-PN spaces, we do always have a topology, but
the family $\{N_q(t)\}$ still generates a certain ``generalized
topology'' as considered by Fr{\'e}chet (see Section
\ref{preliminaries}), where a notion of boundedness is possible
(see Definition \ref{boundedness-GT}), as well as the notion of
$\D$-boundedness is obviously extended (see Definition
\ref{D-boundedness}).

In \cite{LRS3} the authors show that if the PN space is a Serstnev
space which is topologically vectorial, then bounded and
$\D$-bounded subsets coincide. We here extend this result to a
class of pre-PN spaces which we call ``generalized \v Serstnev
spaces'', or more precisely {\it $\phi$-\v Serstnev spaces}, where
$\phi$ is a given map $\R$ to $\R$ satisfying certain properties
(see Theorem \ref{D-bounded-bounded}, and preparatory sections 5
and 7). The well known $\alpha$-simple spaces are indeed $\phi$-\v
Serstnev PN spaces where $\phi(x)=x^{1/\alpha}$ and $\alpha>0$.
This yields a new interpretation of some results in \cite{LRS2};
see Section 6.

Sections 3 and 4 contain some connections between certain PN
spaces with other known structures. For instance, we show that
locally convex topologically vector PN spaces are bornological
(Proposition~\ref{LCTVPN-bornological}). This allows to determine
continuous linear operators between such PN spaces in terms of
bounded subsets (Corollary \ref{operators-alpha-simple}). This
finishes various results and counterexamples in \cite{Iqbal
Jebril} and \cite{LRS3} about linear operators between
$\alpha$-simple spaces. In Theorem~\ref{TVPN-F-normable} we
observe that topologically vector PN spaces are $F$-normable and
paranormable.

We conclude the paper by comparing $\D$-bounded and bounded
subsets in $F$-normable spaces (see Section 9).

\section{Some preliminaries}
\label{preliminaries}

\subsection{PM spaces and PN spaces}
We next recall the definition of PN space given in \cite{alsina1}.
  However, we here do not assume that the triangle
functions involved are continuous. It is convenient to us to
consider also ``triangle functions'' which are non necessarily
associative.

As usual, $\Delta^+$ denotes the set of distance distribution
functions (briefly, a d.d.f.), i.e.~distribution functions with
$F(0)=0$, endowed with the metric topology given by the modified
Levy-Sybley metric $d_S$ (see 4.2 in \cite{SS}). Let $\D^+$
consist of those $F\in \Delta^+$ such that $\lim_{x\to+\infty}
F(x)=1$. Given a real number $a$, $\varepsilon_a$ denotes the
distribution function defined as $\varepsilon_a(x)=0$ if $x\leq a$
and $\varepsilon_a(x)=1$ if $x> a$. Hence, $\R^+$ can be viewed as
a subspace of $\Delta^+$. A triangle function $\tau$ is a map from
$\Delta^+\times \Delta^+$ into $\Delta^+$ which is commutative,
associative, nondecreasing in each variable and has
$\varepsilon_0$ as the identity. If $\tau$ is non associative we
say that it is a {\it non associative triangle function}.

Recall that a {\it probabilistic metric space} (briefly, a PM
space) is a triple $(S,F,\tau)$ where $S$ is a non-empty set, $F$
is a map from $S\times S$ into $\Delta^+$, called the
probabilistic metric, and $\tau$ is an associative triangle
function, such that:
\begin{itemize}
\item[(M1)] $F_{p,q}=\varepsilon_0$ if and only if $p=q$.
\item[(M2)] $F_{p,q}=F_{q,p}$. \item[(M3)]
$F_{p,q}\geq\tau(F_{p,r},F_{r,q})$.
\end{itemize}
When only (M1) and (M2) are required, it is called a {\it
probabilistic semi-metric space} (briefly, PSM space).

A {\it PN space} (respectively, a {\it pre-PN space}) is a
quadruple $(V,\nu,\tau,\tau^*)$ in which $V$ is a vector space
over the field $\R$ of real numbers, the {\it probabilistic norm}
$\nu$ is a mapping from $V$ into $\Delta^+$, $\tau$ and $\tau^*$
are (respectively, neither necessarily commutative nor associative) triangle functions
such that the following conditions are satisfied for all $p$, $q$
in $V$ (we use $\nu_p$ instead of $\nu(p)$):
\begin{itemize}
\item[(N1)] $\nu_p=\varepsilon_0$ if and only if $p=\theta$, where
$\theta$ denotes the null vector in $V$.
\item[(N2)]
$\nu_{-p}=\nu_p$.
\item[(N3)] $\nu_{p+q}\geq\tau(\nu_p,\nu_q)$. \item[(N4)]\,
$\nu_p\leq\tau^*(\nu_{\lambda p},\nu_{(1-\lambda)p})$ for every
$\lambda\in [0,1]$.
\end{itemize}
If, instead of (N1), we only have
$\nu_{\theta}=\varepsilon_0$, then we shall speak of a {\it
(pre-) probabilistic pseudo normed space}, briefly a (pre-) PPN space.

If $(V,\nu, \tau, \tau^*)$ be a PN space (with $\tau$ non
necessarily continuous), then $(V,F, \tau)$ is a probabilistic semi-metric space,
where $F_{p,q}=\nu_{p-q}$.

The following partial order relation is analogous to the corresponding one
for PM spaces (see Section~8.7 of \cite{SS}):
\begin{definition}
\label{better} {\rm A pre-PN space $(V,\nu, \tau_1, \tau_1^*)$ is
{\it better} than another pre-PN space $(V,\nu, \tau_2,
\tau_2^*)$, with the same $V$ and $\nu$, if the following
conditions hold for all $p, q\in V$ and $\lambda\in [0,1]$:
\begin{itemize}
\item $\tau_1(\nu_p, \nu_q)\geq \tau_2(\nu_p,\nu_q)$; \item
${\tau_1}^*(\nu_{\lambda p}, \nu_{(1-\lambda)p})\leq
{\tau_2}^*(\nu_{\lambda p}, \nu_{(1-\lambda) p})$.
\end{itemize}
}
\end{definition}
We do not know if every PN space $(V,\nu, \tau, \tau^*)$ admits a
{\it best-possible} PN structure, in the sense that is better than
any other PN space $(V,\nu, \tau', {\tau^*}')$. It would be
interesting to study this problem for Menger PN spaces (cf.
Section~8.7 in \cite{SS}).

The definition of $\D$-boundedness, which is merely probabilistic
and the same as for PM spaces, was introduced in \cite{LRS1}. It
obviously extends to pre-PN spaces:
\begin{definition}
\label{D-boundedness}
{\rm A subset $A\subseteq V$ in a pre-PN space $(V,\nu,\tau,\tau^*)$ is
called {\it $\mathcal{D}$--bounded} if $\lim_{x\to \infty}R_A(x)=1$, where
$R_A$ is the probabilistic radius  of $A$ given by
\begin{equation}
R_A(x):=l^-\inf\{\nu_p(x):p \in A\}.
\end{equation}
}
\end{definition}

\subsection{Examples of PN spaces}
Recall that a map $T:[0,1]\times [0,1]\to [0,1]$ is a {\it
$t$-norm} if it is commutative, associative, nondecreasing in each
variable, and has $1$ as identity. Then, $\tau_T$ is defined as
$\tau_T(F,G)(x):=\sup\{T(F(s), G(t)): s+t=x\}$, and
$\tau_{T^*}(F,G)(x):=\inf \{T^*(F(s), G(t)): s+t=x\}$, where $T^*$
is the dual $t$-conorm given by $T^*(x,y):= 1-T(1-x, 1-y)$. Notice that if $T$ is
left-continuous then $\tau_T$ is a triangle function
\cite[p.~100]{SS}, although it is not necessary. For instance, if $Z$ is
the minimum $t$-norm, defined as $Z(x,1)=Z(1,x)=x$ and
$Z(x,y)=0$, elsewhere, then $\tau_T$ is a triangle function.

A {\it Menger PM space} under a $t$-norm $T$ is a PM space of the form $(V,
\nu, \tau_T)$. A {\it Menger PN space} (respectively, Menger
pre-PN) under $T$ is a PN  space (respectively, pre-PN
space) of the form $(V, \nu, \tau_T, \tau_{T^*})$.

 A {\it \v Serstnev (pre-) PN space} is a (pre-) PN space
$(V,\nu,\tau,\tau^*)$ where $\nu$ satisfies the following \v
Serstnev condition:
\begin{itemize}
\item[(\v {S})] \quad  $\nu_{\lambda p}(x)=
\nu_p\left(\frac{x}{\mid \lambda\mid }\right)$, for all $x\in
\R^+$, $p\in V$ and $\lambda \in \R\setminus \{0\}$.
\end{itemize}
It turns out that (\v S) is equivalent to have (N2) and
\begin{equation}
\label{characterization-serstnev} \nu_p=\tau_M(\nu_{\lambda
p},\nu_{(1-\lambda)p}),
\end{equation}
for all $p\in V$ and $\lambda\in [0,1]$ (see
\cite[Theorem~1]{alsina1}), where $M$ is the $t$-norm defined as
$M(x,y)=\min\{x,y\}$. Therefore, condition (N4) is satisfied for
every $\tau^*$ such that $\tau_M\leq \tau^*$. In the sense of
Definition \ref{better}, if $(V,\nu, \tau,\tau^*)$ is a \v
Serstnev PN space then $(V,\nu, \tau, \tau_M)$ is a ``better''
structure than $(V,\nu, \tau,\tau^*)$.

\begin{example}
\label{first-example} {\rm Every normed space $(V,\norm{\cdot})$ yields
a \v Serstnev, Menger space $(V,\nu, \tau_M, \tau_M)$, where
$\nu_p=\varepsilon_{\parallel p \parallel}$. (Recall that
$\tau_{M^*}=\tau_M$).}
\end{example}
More generally, when the image of $\nu$ lies in $\R^+\subset
\Delta^+$, i.e $\nu_p=\varepsilon_{g(p)}$ for some function
$g:V\to \R^+$, we obtain in Section 4 the link with $F$-norms.

\begin{example}
\label{def-alpha-simple} {\rm Let $(V,\,\norm \quad)$ be a normed
space, $G\in \Delta^+$ be different from $\varepsilon_0$ and
$\varepsilon_{\infty}$, and $\alpha\geq 0$. Define $\nu: V
\rightarrow \Delta^+$ by $\nu_{\theta}=\varepsilon_0$ and for
$p\neq \theta$:
\begin{equation}
\label{alpha-simple} \nu_p(t):= G\left(
\frac{t}{{\norm{p}^{\alpha}}}\right).
\end{equation}
The triple $(V,G; \alpha)$ is called the {\it $\alpha$--simple
space} generated by $(V, \,\norm{\cdot})$ and~$G$. For $\alpha=1$
we have that $(V, \nu, \tau_M, \tau_M)$ is a Menger PN space, that
is \v Serstnev (see Theorem~2.1 of \cite{LRS2}). If
$G=\varepsilon_1$ we recover Example \ref{first-example}. For
$\alpha\neq 1$ we do not have a Menger PN space in general, but a
PN space with $\tau=\tau^*=\tau_{M,L}$, for some $L$ (see
Proposition~\ref{alpha-simple-PN-spaces}). }
\end{example}

\subsection{Fuzzy normed spaces}
The class of PN spaces has some connection with the class of fuzzy
normed spaces.  We want here establish this connection without
giving many details. This yields a source of examples going in
both directions; see a similar connection for PM spaces and fuzzy
metric spaces in \cite[Section 4]{Olga}. The first definition of
fuzzy norm was given by Katsaras \cite{Kat}, and later extended by
Felbin \cite{Fel}. However, as far as we know, there is only the
article by Wu and Ma \cite{WuMa} relating fuzzy norms and
probabilistic norms as defined in their origin. Recall from
\cite{Fel} that a fuzzy normed space is a quadruple $(V,\norm{\cdot},
L,R)$ where $V$ is a real vector space, $\norm{\cdot}$ is a function
from $V$ to the set of fuzzy numbers, $L$ and $R$ are a continuous
$t$-norm and a $t$-conorm satisfying certain properties.

\begin{proposition}
\label{fuzzy-prob} If $(V,\norm{\cdot}, L,R)$ is a fuzzy normed space,
with $L$ and $R$ continuous $t$-norm and $t$-conorm, respectively,
then $(V, \nu, \tau_{R^*}, \tau_M)$ is a \v Serstnev PN space such
that
$$\nu_{p+q}\leq \tau_{L^*}(\nu_p,\nu_q),$$ for all $p$ and $q$ in $S$,
where $\nu_p$ is the distribution function associated to the fuzzy
number $\norm{p}$.
\end{proposition}

Conversely we have:
\begin{proposition}
\label{prob-fuzzy} If $(V,\nu, \tau_T, \tau_M)$ is a \v Serstnev
PN space and furthermore,
$$\nu_{p+q}\leq \tau_{L^*}(\nu_p,\nu_q),$$ for all $p$ and $q$ in $S$,
then $(V, \norm{\cdot}, L, T^*)$ is a fuzzy normed space.
\end{proposition}

\subsection{The generalized strong topology}
Recall from \cite[Section 12]{SS} (see also \cite{Hoh})
that every PSM
space $(S,F, \tau)$ endows a generalized topology
of type $V_D$ (in the sense of Fr\'{e}chet),
which is Fr{\'e}chet-separated and first-numerable.
It is called the {\it generalized strong topology.}
The associated {\it strong neighborhood system} is given by
$\N=\bigcup_{p\in S} \N_p$, where $\N_p=\{N_p(t): t>0\}$ and
$$N_p(t):=\{q\in S: F_{p,q}(t)> 1-t\}.$$
A countable base at $p$ is given by $\{N_p(1/n): n\in \Na\}$. If
we define $\delta(p,q):= d_S(F_{p,q},\varepsilon_0)$, then
$\delta$ is a semi-metric on $S$, and the neighborhood $N_p(t)$ is precisely
the open ball $\{q: d_S(F_{p,q}, \varepsilon_0)<t\}$.

If $(S,F,\tau)$ is a PM space with $\tau$ continuous, the
generalized strong topology is a genuine topology called the {\it
the strong topology}. Because of (M1) (see subsection 2.1) the
strong topology is Hausdorff-separated. Since it is
first-numerable and uniformable, one has that it is metrizable.


For a pre-PN space $(V,\nu, \tau, \tau^*)$ we have
$N_p(t)=p+N_\theta(t)$, i.e.~the generalized topology is invariant
under translations. The base of $\theta$-neighborhoods
$\{N_\theta(1/n):n\in \Na\}$ determines completely the associated
generalized topology. This is also Fr{\'e}chet-separated, and
countably generated by radial and circled $\theta$-neighborhoods.
In fact a converse result also holds; see \cite{LaR} for more
details.

According to this setting, we give the following definition:
\begin{definition}
\label{boundedness-GT}
{\rm A subset $A\subset V$ in a pre-PN space $(V,\nu,\tau,\tau^*)$ is
{\it bounded\,} if for every integer $m\geq 1$, there is a finite set $A_1\subseteq A$
and a natural number $k\geq 1$ such that
\begin{equation}
 A \subseteq \bigcup_{p\in A_1} (p +
N_\theta(1/m)^{[k]}).
\end{equation}
where here $B^{[k]}=B + \stackrel{(k}{\cdots} + B$.
}
\end{definition}

\subsection{TV groups and TV spaces}
Recall that a vector space endowed with a topology, is a {\it topological vector
space} (briefly, a TV space) if both the addition $ +: V\times
V\to V$ and multiplication by scalars $\eta: \R\times V\to V$ are
continuous. If only the addition is assumed to be continuous then
$V$ is a {\it topological group}; if furthermore $\eta$ is
continuous at the second place, then it is called a {\it
topological vector group} (briefly, a TV group). In~\cite{alsina2}
the authors showed that PN spaces with $\tau$ continuous, are
topological vector groups. We quote this result for further
reference.

\begin{theorem}[\cite{alsina2}]
\label{TVS} A PN space $(V,\nu,\tau,\tau^*)$, with $\tau$
continuous, is a TV space if and only if the map $\eta$ is
continuous at the first place (i.e.~for every fixed $p\in V$,
$\lambda_n p\to 0$ whenever $\lambda_n\to 0$).
 \qed
\end{theorem}
The following conditions to have a TV space are sufficient (see
Theorem~4 and remarks after Theorem~5 in \cite{alsina2}):
$\nu_p\neq \epsilon_\infty$, for all $p\in V$, the subset $\nu(V)$
is closed in $(\Delta^+, d_L)$, $\tau^*$ is continuous, and
$\tau^*$ Archimedean on $\nu(V)$.

For \v Serstnev PN spaces Theorem~\ref{TVS} yields the following
characterization:
\begin{theorem}[\cite{LRS5}]
\label{TV-serstnev}  A \v Serstnev PN space $\VPN$ is a TV space
if and only if $\nu(V)\subseteq \D^+$.\qed
\end{theorem}


If a PN space $(V,\nu,\tau,\tau^*)$ is a TV space
then then a subset is bounded (in the sense of \ref{boundedness-GT}
if and only if for every integer $m\geq 1$,
there is a natural number $k\geq 1$ such that
\begin{equation}
\label{boundedness-TV} A \subseteq k N_\theta(1/m).
\end{equation}
This is also equivalent to being ``topologically bounded'' (as defined in
\cite{alsina1}), that is, for every sequence $(\alpha_n)\subset
\R$ with $\lim_{n} \alpha_n=0$, and for every sequence
$(p_n)\subset A$, then $\lim \alpha_np_n=\theta$ in the strong
topology.

\section{Normable and bornological PN spaces}
Normability of PN spaces has been recently studied in \cite{LRS5}.
The following criterium establishes when a PN space is normable
(see \cite[p.~41]{Schafer}):

\begin{theorem}
(Kolmogorov) \label{Normable} A TV PN space $(V,\nu, \tau,
\tau^*)$ is normable if, and only if, there exists a bounded and
convex $\theta$-neighborhood.
\qed
\end{theorem}
By a Prochaska's result adapted to the theory of PN spaces (see
\cite{LRS5}) we have that all \v Serstnev and Menger PN spaces
$(V,\nu, \tau_M, \tau_M)$ are locally convex.

\noindent
\begin{example}[\cite{LRS5}]
{\rm Let $(V,\norm{\cdot}; G)$ be a simple space and $(V,\nu,\tau_M,\tau_M)$
the associated \v Serstnev and Menger PN space of Example~\ref{def-alpha-simple}.
If $G \in \D^+$ and strictly increasing then the strong topology
is $\norm{\cdot}$-normable.  }
\end{example}

\begin{remark}
{\rm In the previous example, if $G\not\in \D^+$ then the associated strong topology is discrete,
therefore it is not a TV space (it is just a TV group).}
\end{remark}

A locally convex TV space $E$ is \emph{bornological} if every
circled, convex subset $A\subset E$ that absorbs every bounded set
in $E$ is a neighborhood of $\theta$. It is known that metrizable and locally convex
topological vector spaces are bornological (see \cite[II 8.1]{Schafer}).
Bornological spaces are
inductive limits of normable spaces (\cite[II 8.4]{Schafer}).

\begin{proposition}
\label{LCTVPN-bornological} Every PN spaces $(V,\nu, \tau,
\tau^*)$ that is a locally convex TV space is bornological.\qed
\end{proposition}

In Proposition \ref{alpha-simple-PN-spaces} we will see that
$\alpha$-simple spaces are PN spaces.

\begin{example}
\label{alpha-simple-bornological}
{\rm Let $L(x,y)=(x^{1/\alpha}+y^{1/\alpha})^\alpha$.
Then the $\alpha$-simple PN space $(V, \nu, \tau_{M,L}, \tau_{M,L})$ where $\nu_p(x)=
G(x/\norm{p}^\alpha)$, with $G\in \D^+$ is bornological.}
\end{example}

A linear operator $T:V_1\to V_2$ is called {\it bounded}
if it transforms bounded subsets of $V_1$ into bounded subsets of
$V_2$ (see e.g. \cite[p.~63]{Edwards}). Obviously, continuous
linear operators are bounded, but not conversely. However, if the
source space is bornological and the target is a locally convex TV
space then the converse holds (see e.g. \cite[p.~477]{Edwards}).
In particular, we have:

\begin{theorem} A linear operator between two locally convex TV PN spaces
is continuous if, and only if, it is bounded.\qed
\end{theorem}
Example~3.5 in \cite{LRS3} gives a bounded linear operator from a
non bornological (non locally convex) PN space which is not
continuous.

\begin{corollary}
\label{operators-alpha-simple} Let $G$ and $G'$ be in $\D^+$. Let
$(V,G,\alpha)$ and $(V,G',\alpha')$ be two $\alpha$-simple spaces
regarded as PN spaces. A linear operator $T: (V,G,\alpha) \to
(V,G',\alpha')$ is continuous if and only if $T$ is bounded.
\end{corollary}
This corollary closes the results in \cite{Iqbal Jebril} and
\cite[Section~3]{LRS3}.


\section{$F$-normable and paranormable PN spaces}
 Recall from \cite{Schafer} and
\cite[Section~4]{Wil} that an {\it $F$-norm} on a vector space $V$
is a map $g: V\to \R^+$ such that
\begin{itemize}
\item[(i)] $g(p)=0$ if and only if $p=\theta$.
\item[(ii)] $g(\lambda p)\leq g(p)$ if $\abs{\lambda}\leq 1$.
\item[(iii)] $g(p+q)\leq g(p)+g(q)$.
\end{itemize}
The pair $(V,g)$ is called an {\it $F$-normed space}. It is a TV
group with respect to the metric $d(p,q)=g(p-q)$, but in general
it is not a TV space. $F$-normed spaces which are TV spaces are
called {\it paranormed spaces} (see \cite[Section~4]{Wil}).

\begin{example} {\rm Let $V$ be the vector space of all continuous
functions $p:\R\to \R$, $g(p):= \sup_{t\in \R} \frac{ \mid
p(x)\mid}{a+ \mid p(x) \mid}$, with $a>0$. Then $g$ is an $F$-norm
but not a paranorm (see \cite[Exercise~12(b), p.~35]{Schafer}). }
\end{example}

Of course, different $F$-norms may induce the same
metric-topology. For instance, if $(V,\norm{\cdot})$ is a normed space
then $g(p)=\norm{p}^{\alpha}$, or $g(p)=\frac{\parallel p
\parallel}{ \alpha+ \parallel p \parallel}$, where $\alpha>0$,
are $F$-norms which induce the same topology as $\norm{\cdot}$.
Observe that  every $F$-normed (respectively, paranormed) space
$(V,g)$ is homeomorphic to an $F$-normed (respectively, paranormed) space $(V,g')$ with
$g'(V)<1$. Indeed, if $g$ is an $F$-norm, then
$g'(p)=g(p)/(1+g(p))$ is an $F$-norm equivalent to $g$.

The above condition (ii) implies $\norm{-p}=\norm{p}$. This
observation and the fact that
$\tau_M(\varepsilon_a,\varepsilon_b)=\varepsilon_{a+b}$ yield
easily the following correspondence between $F$-norms and certain
PN spaces.

\begin{proposition}
\label{F-norms-PN-spaces} Let $g: V\to \R^+$ be any map and define
$\nu$ by $\nu_p:=\varepsilon_{g(p)}$. Then $(V,g)$ is an
$F$-normed space if, and only if, $(V,\nu,\tau_M, {\bf M})$ is a
PN space, where ${\bf M}$ is defined as ${\bf
M}(F,G)(x)=M(F(x),G(x))$. \qed
\end{proposition}
Notice that ${\bf M}$ is the maximal triangle function, so
$(V,\varepsilon_g,\tau_M,{\bf M})$ could not be the best PN
structure for a given $F$-norm $g$. Indeed, if $g$ is a norm we
can replace ${\bf M}$ by $\tau_M$.

\begin{proposition}
\label{F-norms-PN} Let $g:V\to \R^+$ be any map and define $\nu$
by $\nu_p=\varepsilon_{g(p)}$. Let $\tau$ and $\tau^*$ be two
triangle functions.
\begin{enumerate}
\item
If $\tau(\varepsilon_a,\varepsilon_b)\geq \varepsilon_{a+b}$, for
all\, $a$, $b\in \R^+$, and $(V,\nu,\tau, \tau^*)$ is a PN
space, then $g$ is an $F$-norm.

\item If $\tau(\varepsilon_a,\varepsilon_b)\leq
\varepsilon_{a+b}$, for all\, $a$, $b\in \R^+$, and $g$ is an
$F$-norm, then $(V,\nu, \tau, \tau^*)$ is a PN space if and
only if (N4) holds.


\item
If $\tau(\varepsilon_a,\varepsilon_b)\leq \varepsilon_{a+b}$, for
all\, $a$, $b\in \R^+$, then $g$ is a norm if and only if $(V,\nu,
\tau, \tau^*)$ is a \v Serstnev PN space. \qed
\end{enumerate}
\end{proposition}

\begin{proposition}
Suppose that $(V, \varepsilon_g, \tau,
\tau^*)$ is a PN space, with $g$ an $F$-norm on $V$. Then $\eta$ is continuous at the first place. In this
case, the strong topology is equivalent to the metric-topology
induced by~$g$.
\end{proposition}
\begin{proof}
It is easy to check that the strong neighborhood $N_\theta(t)$
coincides with the open ball $\{p: g(p)<t\}$.
\end{proof}

Conversely, we have the following theorem for TV PN spaces:

\begin{theorem}
\label{TVPN-F-normable} Let $(V,\nu,\tau,\tau^*)$ be a metrizable PN space that is a TV space,
then it is paranormable.
\end{theorem}
\begin{proof}
Theorem~6.1 in \cite[p.~28]{Schafer} implies that metrizable TV
spaces are $F$-normable, therefore by Theorem~\ref{TVS} it is paranormable.
\end{proof}




\section{$\phi$-transforms on PN spaces}
\label{phi-transforms}

Following \cite{alsina0}, for $0<b\leq \infty$, let $M_b$ be the
set of {\it $m$-transforms} which consists on all continuous and
strictly increasing functions from $[0,b]$ onto $[0,\infty]$. More
generally, let $\M$ be the set of non decreasing left-continuous
functions $\phi:[0,\infty]\to [0,\infty]$ with $\phi(0)=0$,
$\phi(\infty)=\infty$ and $\phi(x)>0$, for $x>0$. Then,
$M_b\subset \widetilde M$ once $m$ is extended to $[0,\infty]$ by
$m(x)=\infty$ for all $x\geq b$.  Notice that a function $\phi\in
\M$ is bijective if and only if $\phi\in M_\infty$.

Sometimes, the probabilistic norms $\nu$ and $\nu'$ of two given
(pre-) PN spaces satisfy $\nu'=\nu\phi$ for some $\phi\in \M$, non
necessarily bijective. Let $\phi\h$ be the (unique) quasi-inverse
of $\phi$ which is left-continuous. Recall from \cite[p. 49]{SS}
that $\phi\h$ is defined by $\phi\h(0)=0$, $\phi\h(\infty)=\infty$
and $\phi\h(t)=\sup \{ u : \phi(u)<t \}$, for all $0<t<\infty$. It
follows that $\phi\h(\phi(x))\leq x$ and $\phi(\phi\h(y))\leq y$
for all $x$ and $y$.

One has the following (which generalizes Theorem~4 in \cite{alsina0}):
\begin{theorem}
\label{PPM-phi} Let $(S,F)$ be a PPM space, and $F'=F\phi$
with $\phi\in \M$. Then, $(S,F')$ is a PPM space.
Moreover, the  generalized strong topology induced
by $F$ is finer than the one induced by $F'$. If
$\phi\h(y)>0$, for $y>0$, then they coincide.
\end{theorem}
\begin{proof}
We have that $\phi(x)>0$, for all $x>0$. Hence, for each $m\in\Na$
there is an $n\in \Na$, with $n\geq m$ such that $\phi(1/m)>1/n$.
Thus, for every $p, q\in S$ satisfying $F_{p,q}(1/n)>1-n$, we have
$$F_{p,q}(1/m)=F_{p,q} (\phi (1/m))\geq F_{p,q}(1/n) >1-1/n\geq 1-1/m,$$
i.e.~every strong neighborhood $N'_p(1/m)$ with respect to $F'$
contains a strong neighborhood $N_p(1/n)$ with respect to $F$.
\end{proof}

The following consequences are straightforward:


\begin{corollary}
\label{phi-TV-D-B}
 Let $V_1=(V,\nu,\tau,\tau^*)$ and $V_2=(V, \nu',\tau',(\tau^*)')$ be
two pre-PN spaces with the same base vector space and suppose that
$\nu'=\nu\phi$, for some $\phi\in \M$. Then the following hold:
\begin{enumerate}
\item[(i)] If the scalar multiplication $\eta:\R\times V \to V$ is continuous at the first place
with respect to $\nu$, then it is so with respect to $\nu'$. In
particular, if $\tau$ and $\tau'$ are continuous, and $V_1$ is a
TV PN space, then so is $V_2$.
\item[(ii)] If $\lim_{x\to \infty}\phi(x)=\infty$ and $A$ is a $\D$-bounded in $V_1$ then it so in $V_2$.
\item[(iii)] If $A$ is bounded in $V_1$ then it is so in $V_2$.\qed
\end{enumerate}
\end{corollary}

If $(V,\nu,\tau,\tau^*)$ is a given pre-PN space and $\phi\in \M$,
we can consider the composite $\nu':=\nu \phi$ from $V$ into
$\Delta^+$. By Theorem~\ref{PPM-phi} $\nu'$ satisfies (N1) and (N2). We can
consider the quadruple $(V,\nu \phi,\tau^\phi, (\tau^*)^\phi)$,
where $\tau^\phi$ is given by
\begin{equation}
\label{tau-transform} \tau^\phi(F,G)(x):=\tau(F \phi\h, G
\phi\h)\phi (x),
\end{equation}
and ${\tau^*}^\phi$ is defined in a similar way. The quadruple
$(V,\nu \phi,\tau^\phi, (\tau^*)^\phi)$ is called the {\it
$\phi$-transform} of $(V,\nu,\tau,\tau^*)$.

\begin{proposition}
Let $(V,\nu,\tau,\tau^*)$ be a pre-PN space. If $\phi\in \M$ then
the $\phi$-transform $(V,\nu \phi,\tau^\phi, (\tau^*)^\phi)$ is a
pre-PN space. \qed
\end{proposition}

\begin{remark}
{\rm If $\phi\not\in M_{\infty}$, then associativity of
$\tau^\phi$ and ${\tau^*}^\phi$ might fail. But, if $\phi\in
M_\infty$ then $\tau^\phi$ and ${\tau^*}^\phi$ are (associative)
triangle functions. Hence, in this case the $\phi$-transform of a
PN space is a PN space. Notice also that the $\phi^{-1}$-transform
of $(V, \nu\phi, \tau^\phi, (\tau^*)^\phi)$ is the space
$(V,\nu,\tau,\tau^*)$. } \end{remark}

As in \cite[7.1.7]{SS} let $\Law$ be the set of all binary
operations $L$ on $[0,+\infty]$ which are surjective, non
decreasing in each place and continuous on
$[0,+\infty]\times[0,+\infty]$, except possibly at the points
$(0,+\infty)$ and $(+\infty, 0)$. If $\phi\in M_{\infty}$ and
we define $L(x,y)=\phi^{-1}(\phi(x)+\phi(y))$, then $L\in
\Law$. Given a continuous $t$-norm $T$, one can consider the
triangle functions $\tau_{T,L}$ and $\tau_{T^*,L}$ which are
defined in \cite[7.2]{SS}. An easy calculation yields the
following result:

\begin{theorem} \label{phi-Menger} Let $(V, \nu, \tau_T, \tau_{T^*})$
be a Menger PN space under some continuous $t$-norm $T$, and
$\phi\in M_\infty$. Then, the PN space $(V, \nu\phi, \tau_{T,L}, \tau_{T^*, L})$ is
the $\phi$-transform of $(V, \nu, \tau_T, \tau_{T^*})$.
\end{theorem}
Notice that this is a Menger space under $T$ if, and only if,
$\phi(x)=kx$ for some constant $k\in \R\setminus{0}$ (cf.
\cite[Section~6]{L6}).

\section{$\alpha$-simple PN spaces} \label{alpha-simple-spaces}
As we have seen in Example \ref{def-alpha-simple}, the way to
produce a Menger PN space under $M$ from a simple space
$(V,\norm{\cdot}, G)$ does not need any assumption on the
distribution function $G$.
However, in the case of $\alpha$-simple spaces, some restrictions
on $G$ are required in order to obtain the structure of Menger PN
space under a certain $t$-norm $T_G$ (see Section 3 in
\cite{LRS2}). In this section we give a new proof of
Theorem~3.1, part (a) of \cite{LRS2}, by using the following:

\begin{proposition}
\label{alpha-simple-PN-spaces} If $(V,G,\alpha)$  is an
$\alpha$-simple space, and $\nu_p(t)=G(t/\norm{p}^\alpha)$, then
$(V, \nu, \tau_{M, L}, \tau_{M,L})$ is a PN space, where
$L\in\mathcal{L}$ and $L(x,y)=(x^{1/\alpha}+y^{1/\alpha})^\alpha$.
\qed
\end{proposition}
\begin{proof}
This is a particular case of Theorem~\ref{phi-Menger}, with $\phi(x)=x^{1/\alpha}$.
\end{proof}

Now, suppose that $G\in \Delta^+$ is strictly increasing. Consider
the $t$-norm $T_G$ defined as follows:
$$
T_G(x,y):= G\left(\left\{ [G^{-1}(x)]^{1/(1-\alpha)}+
[G^{-1}(y)]^{1/(1-\alpha)}\right\}^{(1-\alpha)}\right).
$$

\begin{corollary}[\cite{LRS2}]Let $(V, G,\alpha)$ be an $\alpha$-simple space, where $G$ is an
strictly increasing continuous distribution function, and
$\alpha>1$. Then $(V,\nu,\tau_{T_G}, \tau_{{T_G}^*})$ is a Menger
PN space under $T_G$.
\end{corollary}
\begin{proof}
Let $\tau=\tau_{T_G}$ in the above proposition. We have to see
that $(V,\nu, \tau_{M,L}, \tau_{M,L})$ is better than
$(V,\nu,\tau, \tau^*)$ in the sense of Definition \ref{better}.
For that, we have to show that
$\tau_{M,\,L}(\nu_p,\nu_q) \geq \tau(\nu_p,\nu_q)$ and
$\tau_{M,\,L}(\nu_{\lambda p}, \nu_{(1-\lambda)p})
\leq \tau^*(\nu_{\lambda p}, \nu_{(1-\lambda)p})$,
for all $p, q \in V$ and $\lambda \in (0,1)$.
$$
\begin{array}{lll}
\tau(\nu_p,\nu_q)(x)=\sup_{r+s=x} \{T_G(\nu_p(r),\nu_q(s))\} =&

\\ \sup_{r+s=x}\left\{G \left( \left[G^{-1}(G\left(\frac{r}{\parallel p \parallel^\alpha
 }\right)\right]^{1/(1-\alpha)}+
 \left[G^{-1}\left(G\left(\frac{s}{\parallel q\parallel^\alpha}
\right)\right)\right]^{1/(1-\alpha)} \right) ^{(1-\alpha)}
\right\}=&\\

\sup_{r+s=x}\left\{G\left(\left(\frac{r}{\parallel p
\parallel^\alpha
 }\right)^{1/(1-\alpha)}+
  \left(\frac{s}{\parallel q\parallel^\alpha}
\right)^{1/(1-\alpha)}\right)^{(1-\alpha)} \right\}.
\end{array}
$$
On the other hand
$$
\begin{array}{lll}

\tau_{M,\,L}(\nu_p,\nu_q)(x)=\sup_{L(u,
v)=x}\left\{M(\nu_p(u),\nu_q(v))\right\}
=&\\
\sup_{u^{1/\alpha}+ v^{1/\alpha}=x^{1/\alpha}}\left\{M\left(G
\left(\frac{u}{\parallel p
\parallel^\alpha}\right),G\left(\frac{v}{\parallel q
\parallel^\alpha}\right)\right)\right\}=&\\
G\left(\frac{x}{(\parallel p
\parallel+\parallel q
\parallel)^\alpha}\right).

\end{array}
$$
Now, we use one of the known H{\"o}lder's inequalities
$$ (a+b)^{1-\alpha}\leq
\lambda^\alpha a^{1-\alpha} + (1-\lambda)^\alpha b^{1-\alpha},$$
which holds for $\alpha>1$, $\lambda\in(0,1)$ and $a, b\in (0,+\infty)$.
By setting $$a:=\left(\frac{r}{\parallel p
\parallel^\alpha
 }\right)^{1/(1-\alpha)}\qquad b:=\left(\frac{s}{\parallel q\parallel^\alpha}
\right)^{1/(1-\alpha)}\qquad \lambda:=\frac{\parallel p
\parallel^\alpha}{(\parallel p
\parallel  + \parallel q
\parallel)^\alpha}$$
it follows
$$
\left(
\frac{r}{\parallel p \parallel^\alpha}
\right)^{1/(1-\alpha)}
+
\left(\frac{s}{\parallel q \parallel^\alpha}
\right)^{1/(1-\alpha)}
\leq
\left( \frac{r+s}{\parallel p \parallel  + \parallel q \parallel)^\alpha}
\right)^{1/(1-\alpha)}
$$
After applying $G$ in both sides, we obtain one of the desired inequailties
$\tau_{M,\, L}(\nu_p,\nu_q) \geq \tau(\nu_p,\nu_q)$.
The other inequality follows analogously.
\end{proof}

A similar result can be shown for $\alpha<1$ by choosing a
$t$--\emph{norm} $T$ as in Theorem~3.2, part (a) of \cite{LRS2}.

\section{$\phi$-transforms on Serstnev spaces}
\label{phi-Serstnev}

If $(V,\nu, \tau,\tau^*)$ is a \v Serstnev pre-PN space and
$\nu':=\nu \phi$ for some bijective function $\phi\in M_{\infty}$
(see the Section~5).  Then $\nu'$ will satisfy $\nu'_{\lambda
p}(x)= \nu'_p\left(\phi^{-1}\left(\frac{\phi(x)} {\mid \lambda\mid
}\right)\right)$, for all $x\in \R^+$, $p\in V$ and $\lambda \in
\R\setminus \{0\}$. This motivates the following definition for
$\phi$ non necessarily bijective.

\begin{definition}
{\rm We say that a quadruple $(V,\nu,\tau,\tau^*)$ satisfies the
$\phi$-\v Serstnev condition if:
\begin{itemize}
\item[($\phi$-\v {S})] $\quad \nu_{\lambda p}(x)=
\nu_p\left(\phi\h\left(\frac{\phi(x)} {\mid \lambda\mid
}\right)\right)$, for all $x\in \R^+$, $p\in V$ and $\lambda \in
\R\setminus \{0\}$.
\end{itemize}
A pre-PN space $(V,\nu,\tau,\tau^*)$ which satisfies the $\phi$-\v
Serstnev condition is called a {\it $\phi$-\v Serstnev pre-PN
space}.}
\end{definition}

\begin{example}
{\rm If $\phi(x)=x^{1/\alpha}$, for a fixed positive real number
$\alpha$, then condition ($\phi$-\v S) takes the form:
\begin{itemize}
\item[($\alpha$-\v S)] $\quad \nu_{\lambda p}(x)=
\nu_p\left(\frac{x}{\mid \lambda\mid^\alpha} \right)$, for all
$x\in \R^+$, $p\in V$ and $\lambda \in \R\setminus \{0\}$.
\end{itemize}
Pre-PN spaces satisfying ($\alpha$-\v S) are called {\it
$\alpha$-\v Serstnev}. (Note $1$-\v Serstnev is just \v Serstnev.)
We will see in Proposition~\ref{alpha-simple-PN-spaces} that
$\alpha$-simple spaces give rise to $\alpha$-\v Serstnev PN spaces
of the form $(V,\nu, \tau_{M,L}, \tau_{M,L})$ where
$L(x,y):=(x^{1/\alpha}+ y^{1/\alpha})^\alpha$. Thus,
$\alpha$-simple pre-PN spaces can be viewed as $\phi$-transforms
of PN simple spaces.}
\end{example}
More generally, the $\phi$-transform of a \v Serstnev PN space is
a $\phi$-\v Serstnev pre-PN space, if $\phi$ is bijective. This
yields the following characterization for $\phi$-\v Serstnev
pre-PN spaces.

\begin{proposition}
\label{characterization-phi-serstnev} Let
$L(x,y)=\phi^{-1}(\phi(x)+\phi(y))$ with $\phi\in M_\infty$. Then
($\phi$-\v S) holds if and only if (N2) and also
$$\nu_p=\tau_{M,L}(\nu_{\lambda p}, \nu_{(1-\lambda)p}),$$
for every $p\in V$ and $\lambda \in [0,1]$ are satisfied. In
particular, $\phi$-\v Serstnev spaces admit a better pre-PN
structure of the form $(V,\nu,\tau, \tau_{M,L})$.
\end{proposition}
\begin{proof}
By \cite[Section~7.7]{SS}, taking quasi-inverses we have that
$$\nu_p=\tau_{M,L}(\nu_{\lambda p}, \nu_{(1-\lambda)p})\Longleftrightarrow
\nu_p\h=L(\nu_{\lambda p}\h, \nu_{(1-\lambda) p}\h).$$ By
definition of $L$, this is equivalent to $\phi \nu_p\h=\phi
\nu_{\lambda p}\h + \phi \nu_{(1-\lambda) p}\h$. Taking again
quasi-inverses, we obtain $\nu_p\phi^{-1}= \tau_M(\nu_{\lambda
p}\phi^{-1}, \nu_{(1-\lambda) p}\phi^{-1})$. This condition
together with (N2) is equivalent to the \v Serstnev condition for
$\nu\phi$.
\end{proof}

\noindent For $\alpha$--\v {S}erstnev spaces one also has a
slightly different characterizing formula.

\begin{proposition}
\label{characterization-alpha-serstnev} Let $\alpha\in \R^+$. Then
{\rm ($\alpha$-\v S)} holds if, and only if, (N2) and
\begin{equation}
\label{N4-equation} \nu_{\beta p}= \tau_M(\nu_{\lambda
p},\nu_{(1-\lambda)p}),
\end{equation}
for every $p\in V$ and $\lambda \in [0,1]$ are satisfied, where
$\beta=[\lambda^{\alpha}+(1-\lambda)^{\alpha}]^{1/\alpha}$.
\end{proposition}
\begin{proof} Suppose first that ($\alpha$-\v S) is satisfied.
Then, obviously $\nu_{-p}=\nu_{p}$, hence (N2) holds. As in the
proof of the previous proposition, we have that
(\ref{N4-equation}) holds if, and  only if,
 \begin{equation}
 \label{sumadenus}
 \nu_{\beta p}\h=\nu_{\lambda
 p}\h+\nu_{(1-\lambda)p}\h,
\end{equation}
for all $p\in V$ and all $\lambda \in [0,1]$. Then, because of
 ($\alpha$--\v {S}),
$$\nu_{\lambda p}(x)=\nu_p\left(\frac{x}{\lambda^\alpha}\right)
\Longleftrightarrow \nu_{\lambda
 p}\h= \lambda^{\alpha} \nu_p\h,$$ for every $\lambda\in [0,1]$ and for every
$p\in V$, so that (\ref{N4-equation}) holds easily.

Conversely, suppose that (N2) and (\ref{N4-equation}) hold. Let
$g:\R^+ \rightarrow\R$ be defined as
 $$g(z):=\nu_{zp}\h(t),$$ for fixed $t\in [0,1]$ and $p\in V$.
Then $g$ is a non-decreasing map such that
$$g[(\lambda^\alpha+(1-\lambda)^\alpha)^{1/\alpha}z]=g(\lambda z)+g[(1-\lambda)z].$$
Define now $f:\R^+\to \R$ by $f(z):=g(z^{1/\alpha})$. Then, $f$ is
a nondecreasing function that satisfies $f(x+y)=f(x)+f(y)$, for
all $x$ and $y\in \R^+$. This is the Cauchy equation, therefore by
\cite[Corollary 5]{Aczel} we have $f(x)=f(1)\cdot x$, that is
$g(x^{1/\alpha})=g(1)\cdot x$. By taking $z=x^{1/\alpha}$, we
obtain $g(z)=g(1)z^{\alpha}$, and hence $\nu_{zp}\h(t)= z^{\alpha}
\nu_{p}\h(t)$. This last equality yields ($\alpha$-\v S), as
desired.
\end{proof}

 \begin{corollary}
 \label{new-alpha-serst}
If $(V,\nu,\tau,\tau^*)$ is an $\alpha$--\v {S}erstnev space, then
$(V,\nu,\tau,\tau_M)$ is also an $\alpha$--\v {S}erstnev space.
 \end{corollary}
\begin{proof}
This follows from the inequalities $\nu_p(x) \leq \nu_{\beta p} (x)=\tau_M(\nu_{\lambda p},\nu_{(1-\lambda)p})$.
\end{proof}

\section{Boundedness in $\phi$-\v Serstnev spaces}

An example in \cite{LRS3} shows that
even for normable PN spaces bounded subsets and $D$-bounded
subsets do not coincide. They do coincide on the \v Serstnev PN spaces
that are TV spaces (\cite[Theorem~2.3]{LRS3}). In fact, for \v Serstnev spaces
which are not TV spaces, $\D$-bounded subsets are bounded,
but the converse might fail as we illustrate with Example~\ref{bounded-not-D-bounded}.


We first generalize the following:
\begin{theorem}[\cite{LRS5}]
\label{TV-Serstnev}  A \v Serstnev PN space $\VPN$ is a TV space
if, and only if, $\nu(V)\subseteq \D^+$.\qed
\end{theorem}

\begin{theorem}
\label{TV-phi-Serstnev} Let $\phi\in \M$ such that $\lim_{x\to
\infty} \phi\h(x)=\infty$. Let $(V,\nu,\tau,\tau^*)$ be a
$\phi$-\v Serstnev pre-PN space. Then scalar multiplication
$\eta:\R\times V\to V$ is continuous at the first place if and
only if $\nu(V)\subseteq\mathcal{D}^+$.
\end{theorem}

\begin{proof}
If $\nu$ maps $V$ into $\mathcal{D}^+$,
then, for every $x>0$ and every sequence $\{\alpha_n\}$ converging
to $0$,  one has:
$$\nu_{\alpha_n p}(x)= \nu_p\left(\phi\h\left(\frac{\phi(x)} {\mid \alpha_n\mid
}\right)\right) \longrightarrow 1,$$ as $n$ tends to $+\infty$ (we
use the fact that $\lim_{y\to \infty} \phi\h(y)=\infty$), whence
the assertion.

Conversely, suppose that $\eta$ is continuous at the first place. For
every $n\geq 1$, let $x_n=\phi\h(n\phi(1))$. Then, for all $p\in V$,
$$
\nu_p(x_n)= \nu_p(\phi\h(n\phi(1))=
\nu_p\left(\phi\h\left(\frac{\phi(1)}{1/n}\right)\right) =
\nu_{p/n}(1) \longrightarrow 1.$$ The last term converges to 1
by assumption. Therefore, $\nu_p(x)\to 1$ whenever $x$
tends to infinity, as desired.
\end{proof}

A remarkable result in \cite{LRS3} is Theorem 2.3, where it is
shown that in a \v Serstnev space that is a TV space, a subset is
$\D$-bounded if, and only if, it is bounded or ``topologically
bounded'' (fact that it has been observed in the introduction to be same). We
extend this result to $\phi$-\v Serstnev spaces in the following
theorem with almost the same proof. Notice that the implicit
assumption in \cite{LRS3} that they are TV spaces is not necessary
at all for the first part. The restriction to TV spaces
generalizes a result in \cite{LRS5}.

\begin{theorem}
\label{D-bounded-bounded} Let $\phi\in \M$ such that $\lim_{x\to
\infty} \phi\h(x)=\infty$. Let $(V,\nu, \tau, \tau^*)$ be a
$\phi$--\v Serstnev pre-PN space. Then for a subset $A\subset V$
the following statements are equivalent:
\begin{itemize}

\item[(a)] For every $n\in \Na$, there is a $k\in \Na$ such that $
A\subset k N_{\theta}(1/n). $

\item[(b)] $A$ is $\mathcal{D}$--bounded.
\end{itemize}
These equivalent conditions imply:
\begin{itemize}
\item[(c)] $A$ is bounded.
\end{itemize}
In particular, a subset of a $\phi$-\v Serstnev PN space that is a TV space
is $\D$-bounded if and only if it is bounded.
\end{theorem}
\begin{proof}
Suppose that (a) holds. For every $n\in \Na$, there is a $k\in
\Na$ such that $\nu_{p/k}(1/n)>1-1/n$ for all $p\in A$. Since
$\phi$ is non-decreasing and continuous at infinity, there exists
an $x_0\in\R^+$ such that for all $x\geq x_0$,
$\phi\h(\phi(x)/k)\geq 1/n$. Thus, for every $p\in A$ and $x\geq
x_0$, we obtain
$$
\nu_p(x)=\nu_{k \frac{p}{k}}(x)= \nu_{\frac{p}{k}}\left(\phi\h
\left(\frac{\phi(x)}{k}\right)\right) \geq
\nu_{\frac{p}{k}}\left(\frac{1}{n}\right)>1-\frac{1}{n},
$$
so that, $R_A(x)\geq 1-1/n$, i.e.~$R_A\in \D^+$ as desired.

Conversely, suppose that $A$ is $\D$-bounded. Then, for every
$n\geq 1$ there is an $x_n>0$ such that $R_A(x_n)>1-1/n$. Thus,
$\nu_p(x_n)\geq R_A(x_n)> 1-1/n$, for all $p\in A$. As before,
there exits a $k\in \Na$ such that $\phi\h(k\phi(1/n))\geq x_n$.
Then, for all $p\in A$,
$$\nu_{\frac{p}{k}}\left(\frac{1}{n}\right)=
\nu_{p}\left(\phi\h\left(k\phi\left({\frac{1}{n}}\right)\right)\right)
\geq \nu_{p}(x_n)>1-{\frac{1}{n}},$$ as desired. Finally, (a)
implies (c) because $k N_\theta(1/n)$ is contained in~$N_\theta(1/n)^{[k]}$.
\end{proof}

\begin{example}
\label{bounded-not-D-bounded} {\rm If $(V,\norm{\cdot}, G)$ is a
simple space with $G\not\in \D^+$, then it is a \v Serstnev PN
space and its topology is discrete,
thus not a TV space. In this case, a single set $\{p\}$, with
$p\in V\setminus \{\theta\}$, is bounded but not $\D$-bounded. }
\end{example}
\section{Boundedness in $F$-normable PN spaces}
We include a section treating the following unsolved problem:

\begin{problem} \label{Hohle-problem2} {\rm Determine the
class of all TV PN spaces where $\D$-bounded and bounded subsets
coincide.}
\end{problem}

\begin{proposition}
Such spaces satisfy that $\nu(V)\subseteq \D^+$.
\end{proposition}
\begin{proof} Indeed, for every $p\in V$, the map $\R\to V$, given by
$\lambda\mapsto \lambda p$, is continuous. This implies $\{p\}$
bounded. Hence, $\nu_p\in \D^+$.
\end{proof}

However, the condition $\nu(V)\subseteq \D^+$ is not sufficient to
have the equivalence between boundedness and $\D$-boundedness (see example
below).


\begin{proposition}
Let $(V,g)$ be an $F$-normed space and $(V,\nu,\tau, \tau^*)$ be
any pre-PN space with $\nu_p=\varepsilon_{g(p)}$. Then $A$ is
$\D$-bounded if and only if $g(A)$ is bounded in~$\R^+$.
\end{proposition}
\begin{proof}
Suppose that $A$ is not $\D$-bounded, then $\lim_{x\to \infty}
R_A(x)\not= 1$, hence this limit must be $0$. Hence, for every
$k\geq 1$ there exists a $p_k\in A$ such that
$\varepsilon_{g(p_k)}(k)=0$. This implies $g(p_k)\geq k$ for all
$k\geq 1$, and therefore $g(A)$ is unbounded. The converse can be
proved similarly.
\end{proof}

\begin{proposition}
Let $(V,g)$ be an $F$-normed space. Then:
\begin{enumerate}
\item If $A$ is $\D$-bounded, so is $kA$ for all $k\in \R^+$.
\item Every bounded
subset is $\D$-bounded.
\end{enumerate}
\end{proposition}
\begin{proof} For the first part, suppose that $kA$ is not
$\D$-bounded, for a natural number $k$. Then, there exists a
sequence $(kp_r)\subseteq kA$ with $p_r\in A$ and $g(kp_r)$ converging to
infinity. Since, $g(kp_r)\leq g(p_r)+ g((k-1)p_r)$, we have that
either $g(p_r)$ or $g((k-1)p_r)$ tends to infinity. By induction
we can obtain that $g(p_r)$ tends to infinity.

For the second part, let $A$ be a bounded subset of $V$. Suppose
that $A$ is not $\D$-bounded. Then, there exists a sequence
$(p_r)\subseteq A$ with $g(p_r)$ converging to infinity. Since
$(p_r)$ is bounded, given $n=1$ there exists a $k$ such that
$g(p_r/k)<1$ for all $r\geq 1$. But by part 1, $g(p_r/k)$
converges to infinity, which is a contradiction.
\end{proof}

\begin{example}
{\rm Consider the PN space $(\R, \nu, \tau_M, {\bf M})$ (see
Theorem~5 and Example~4 in \cite{L5}) where
$\nu_p=\varepsilon_{\parallel p
\parallel/(1 + \parallel p\parallel)}$.
Then $\nu_p\geq \varepsilon_1\in \D^+$ for all $p\in \R$. Thus
$\R$ is $\D$-bounded, but of course it is not bounded.}
\end{example}


Another open problem related to problem~\ref{Hohle-problem} is
the following:
\begin{problem}
{\rm Determine the class of all PN spaces $(V, \nu,\tau,\tau^*)$,
with $\tau^*$ Archi\-medean (thus TV spaces), where $\D$-bounded
and bounded subsets coincide.}
\end{problem}

\bigskip
{\bf Acknowledgments} We thank C. Sempi for useful comments, and
B. Schweizer to suggest us Theorem \ref{TVPN-F-normable}. The
research of the first author was supported by grant from the
Ministerio de Ciencia y Tecnolog\'{\i}a (BFM2003-06522) and from
the Junta de Andaluc\'{\i}a (CEC-JA  FQM-197). The second author
is supported by the Junta de Andaluc\'{\i}a (CEC-JA FQM-213).


\end{document}